\documentclass[10pt]{article}

\usepackage{latexsym,amssymb,amsmath,graphicx, amsthm, times}%,bbm}
\usepackage{times}

\newcommand {\CA}{\mathfrak{A}}

\newcommand {\ME}{\mathbb{E}^{x}}

\renewcommand{\P}{\mathbb{P}}

\numberwithin{equation}{section}
\newtheorem{proposition}{Proposition}[section]
\newtheorem{corollary}[proposition]{Corollary}
\newtheorem{remark}{Remark}[section]
\newtheorem{lemma}[proposition]{Lemma}
\newtheorem{example}{Example}[section]
\newtheorem{assump}{Assumption}[section]

\newcommand {\R}{\mathbb{R}}

\newcommand {\F}{\mathcal{F}}
\newcommand {\A}{\mathcal{A}}

\newcommand {\p}{\mathbb{P}}

\newcommand {\E}{\mathbb{E}}

\title{An Analysis of Monotone Follower Problems for Diffusion Processes
}
\author{Erhan Bayraktar \thanks{E. Bayraktar is supported in part by the National Science Foundation, under grant
DMS-0604491.}
\and Masahiko Egami% <-this % stops a space
\thanks{E. Bayraktar and M. Egami are in the Department of Mathematics, University
of Michigan, Ann Arbor, MI 48109, USA, email: \{erhan,
egami\}@umich.edu. } }
\date{}
\begin{document}

\maketitle
\begin{abstract}\noindent

We consider a singular stochastic control problem, which is called
the Monotone Follower Stochastic Control Problem and give
sufficient conditions for the existence and uniqueness of a
local-time type optimal control. To establish this result we use a
methodology that has not been employed to solve singular control
problems. We first confine ourselves to local time strategies.
Then we apply a transformation to the total reward accrued by
reflecting the diffusion at a given boundary
 and
show that it is linear in  its continuation region. Now, the
problem of finding the optimal boundary becomes a non-linear
optimization problem: The slope of the linear function and an
obstacle function need to be simultaneously maximized. The
necessary conditions of optimality come from first order
derivative conditions. We show that under some weak assumptions
these conditions become sufficient. We also show that the local
time strategies are optimal in the class of all monotone
increasing controls.

As a byproduct of our analysis, we give sufficient conditions for
the value function to be $\mathbf{C}^2$ on all its domain.
 We
solve two dividend payment problems to show that our sufficient
conditions are satisfied by the examples considered in the
mainstream literature. We show that our assumptions are satisfied
not only when capital of a company is modeled by a Brownian motion
with drift but also when we change the modeling assumptions and
use a square root process to model the capital.
%We also sketch how
%our approach can be generalized to Bounded Variation Follower
%Stochastic Control Problem and give the characterization of the
%optimal boundaries. We solve two inventory management problems. %In
%the first problem the contents of the inventory fluctuate as a
%Brownian motion with drift and in the second one the contents
%fluctuate as an Ornstein-Uhlenbeck process.
\end{abstract}

\section{Introduction}

We solve a class of singular control problems which are known as
Monotone Follower Stochastic Control Problems (see Karatzas and
Shreve \cite{Karat-Shr1984} for the terminology) for a general
class of diffusion processes. In particular, we give necessary and
sufficient conditions under which the continuation region is
constituted by a single open interval in the state space of the
controlled process. To establish our main result, we first
restrict ourselves to local-time strategies, each of which
corresponds controlling the underlying diffusion by reflecting it
at a particular point. Applying a particular transformation to the
 total reward
accrued by reflecting the diffusion at a given boundary, we show
that the transformed reward is linear in its continuation
region.The slope is a function of the boundary point. In the rest
of the state space, in the region of action, the transformed
reward is equal to an obstacle, which also depends on the value of
the boundary point. This transforms finding the optimal boundary
to a non-linear optimization problem: The slope and the obstacle
have to be maximized simultaneously. We give the necessary
conditions of optimality using the first order derivative
conditions and show that under some weak assumptions these
conditions become sufficient. That is, our methodology of
identifying the
 unique solution of the singular control problem relies
 on a combination of the classical diffusion theory,
  which helps us give a geometric characterization of the value function (the optimal reward), and non-linear programming.
Next, we show that the local time strategies are optimal in the
class of monotone increasing strategies under some certain
assumptions.
% We also show how our
%analysis can be extended to the bounded variation follower problem
%(two-sided singular control problem) and give the characterization
%of the optimal boundaries.

Among the benefits of our analysis are the sufficient conditions
we provide for the value function to be $\mathbf{C}^2$ in the
entire state space. This sheds light on the heuristic
\emph{principle of smooth fit}, which suggests that the value
function is $\mathbf{C}^2$ across the boundary that demarcates the
regions of action and inaction. Our approach should be contrasted
with the ad hoc ordinary Hamilton-Jacobi-Bellman (HJB) approach,
which assumes the principle of smooth fit to construct a solution.
There is no guarantee that a solution could be found and using
that method it is hard to establish sufficient conditions under
which a solution exists. And even if a solution is constructed to
the quasi-variational inequalities, one still has to verify
whether the solution satisfies the assumptions of a
\emph{verification} lemma, i.e., verify the optimality. For
further details of this
approach see e.g. \O ksendal and Sulem \cite{oksendal-book-2}.% On
%the other hand, when one tries to solve the monotone follower
%problem with the ordinary HJB approach, she needs to solve two
%equations for two unknowns, one of which is the optimal boundary
%and the other is one of the undetermined coefficients of the
%homogeneous solution of the HJB equation. This should be
%contrasted with our characterization of the optimal boundary as
%the solution of a single non-linear equation. We stress that we do
%not rely on the assumption that the value function is
%$\mathbf{C}^2$. After computing the value function we find out
%that this is indeed the case under our assumptions. In the bounded
%variation follower problem, the ordinary HJB equation method
%produces four equations for four unknowns whereas our
%characterization gives us two equations for two unknowns.

To illustrate our results we consider the dividend payment problem
for two different scenarios. First, we take the cash-flow of a
company to be a Brownian motion with drift. (This case was
analyzed by Jeanblanc and Shiryaev \cite{JBS1995} using the
ordinary HJB approach.) Second, we take the cash flow of the
company to be a square root process. In this case we show that the
optimal reflection level is strictly less than the mean-reversion
level. In the second example the functions in terms of which the
sufficient conditions are stated are only available in terms some
special functions. Yet, we are able to prove that the sufficient
conditions in our theorems are satisfied by only analyzing the
ordinary differential equation these functions satisfy without
making a reference to their explicit representations. This gives
us a method to check the sufficient conditions for other
diffusions even when explicit representations are not available.
We also extend our results to solve constrained optimization
problems. A dividend payout problem with solvency constraints was
recently solved by Paulsen \cite{paulsen}. In this problem, the
firm is allowed to pay dividends only if the cash flow process is
greater than certain (pre-determined) value. Here, we provide a
simpler proof to Paulsen's result, by a very simple modification
of the proof of Proposition 2.3, which characterization provided
for the optimal reflection barrier.

A similar methodology to ours was used by Dayanik and Karatzas
\cite{DK2003}, to give a general characterization of the value
function of the optimal stopping problem of one dimensional
diffusions. The value function of the optimal stopping problem (up
to a transformation) is characterized as a concave majorant of a
\emph{fixed} obstacle. In the singular control problem we analyze,
the obstacle is not fixed. When we apply the same transformation
to the reward corresponding to the control that is identified by a
fixed boundary, the transformed reward becomes linear in the
region of inaction whose slope depends on this boundary point. On
the region of action the transformed reward is an obstacle and is
a function of the boundary point. Therefore, we maximize the slope
of the linear function and the obstacle simultaneously over all
possible boundary points to obtain the optimal boundary. As a
result, we characterize the optimal boundary first and compute the
value function (the optimal reward) given
 this characterization, whereas \cite{DK2003} characterize
 the value function first and then compute the optimal boundary using this characterization.

Dayanik and Egami \cite{DE2005},  Bayraktar and Egami
\cite{BE2006} (in this work effects of implementation delay are
taken into account) also use similar methodologies. However, the
results we obtained here can not be obtained from the results of
these papers. In these papers, we could not characterize the
optimal control policy completely. The boundary of the region of
action and inaction can be completely characterized only when the
threshold that the state process is taken to, after the
application of the control, is given. Therefore, the
characterization of the optimal boundary that
 we obtain here for the singular control problem
 can not be obtained using a limiting argument (as the fixed cost goes to zero).
Also, the two problems are very different in nature. For example,
the singular control problem is smoother than the impulse control
problem. In the impulse control problem, given a particular
policy, both the slope of the transformed reward in the region of
inaction and the obstacle (the transformed value function in the
region of action) can be determined using the fact that it is
continuous at the boundary.   However, determining the slope of
the transformed reward in the singular control problem is
trickier. To write down the slope of the transformed reward in the
region of inaction as a function of the boundary, we first show
that the transformed reward is $\mathbf{C}^1$. The continuity of
the first derivative is also used to determine the obstacle (the
transformed reward in the region of action) as a function of the
boundary point. On the other hand, the proof of optimality of
local time strategies among a more general class of controls in
the case of singular control problem differs significantly from
the optimality proof of the threshold strategies in the case of
impulse control problem. The latter uses the fact that the value
function (optimal reward) of the impulse control problem can be
approximated by a sequence of optimal stopping problems.

The rest of the paper is organized as follows:  In section 2, we
solve the monotone follower problem for a general diffusion. We
first find the optimal local time strategy. In Section 2.1, we
state the problem, in Section 2.2, we characterize the value
function corresponding to a given boundary and after applying a
particular transformation value function becomes linear in the
region of inaction. In Section 2.3, we characterize the optimal
local time control. We also extend our analysis to solve a
constrained optimization problem. In Section 2.4, we show that the
local time strategies are optimal among all \emph{admissible}
monotone controls. Here, we also point out that under the
assumptions of Proposition~\ref{prop:singular-solution} the value
function is $\mathbf{C}^2$. In Section 2.5, we solve the dividend
payment problem when the capital of a company is assumed to be
either Brownian motion with drift or a square root process.
% we discuss how to extend
%our results to bounded variation control problems and give a
%characterization of the optimal control.
We collect some preliminary results to Section 3, which is our
appendix.

 %Here, we also solve the
%inventory control problem for a Brownian motion with drift and for
%an Ornstein-Uhlenbeck process.

\section{Solution of Monotone Follower Problems}
\subsection{Reflected Diffusions}
 Let $(\Omega,
\F, (\F_t)_{t\geq 0}, \p)$ be a complete probability space with a
standard Brownian motion $W=\{W_t; t\geq 0\}$ and consider the
diffusion process $X^0$ with state pace $\mathcal{I}=[c,
d)\subseteq \mathbb{R}$ and dynamics
\begin{equation}\label{eq:process}
dX^0_t=\mu(X^0_t)dt + \sigma(X^0_t)dW_t
\end{equation}
for some Borel functions $\mu :\mathcal{I}\rightarrow \mathbb{R}$
and $\sigma :\mathcal{I}\rightarrow (0, \infty)$. (We assume that
the functions $\mu$ and $\sigma$ are sufficiently regular so that
(\ref{eq:process}) makes sense.) %Here we take $c$ and $d$ to be a
%natural boundary (the process may start from $c$ but never reaches
%$c$ again) and $d$ to be an absorbing boundary.
We use ``0" as the superscript to indicate that $X^0$ is
uncontrolled. We denote the infinitesimal generator of $X^0$ by
$\A$ and consider the ODE $(\A-\alpha)v(x)=0$. This equation has
two fundamental solutions, $\psi(\cdot)$ and $\varphi(\cdot)$. We
set $\psi(\cdot)$ to be the increasing and $\varphi(\cdot)$ to be
the decreasing solution.  \footnote{In fact, defining $\tau^0_r
\triangleq \inf\{t \geq 0: X^0_t=r\}$, for every $r \in (c,d)$, we
can write
\begin{equation}\label{eq:var-phi-psi}
\psi(x)=\begin{cases}\mathbb{E}^{x}\left[e^{-\alpha
\tau^0_y}\right],& \text{if}\, x \leq y
\\ 1/\mathbb{E}^{y}\left[e^{-\alpha
\tau^0_x}\right], &\text{if}\, x > y
\end{cases},
\quad \varphi(x)=\begin{cases} 1/\mathbb{E}^{y}\left[e^{-\alpha
\tau^0_x}\right],& \text{if}\, x \leq y,
\\\mathbb{E}^{x}\left[e^{-\alpha \tau^0_y}\right],&
\text{if}\, x > y
\end{cases},
\end{equation}
for every $x \in \mathcal{I}$ and an arbitrary but fixed $y \in
\mathcal{I}$ (see It\^{o} and McKean \cite{ito-mc}).} We will take
$c$ to be absorbing and $d$ to be natural, and therefore
$\psi(d-)=\infty$, $\varphi(d-)=0$ since $X^0$ never reaches $d$.
First, we define an increasing function
\begin{equation}\label{eq:F}
F(x)\triangleq\frac{\psi(x)}{\varphi(x)}.
\end{equation}
Next, we define concavity of a function with respect $F$ as
follows: A real valued function $u$ is called \emph{$F$-concave}
on $(c, d)$ if, for every $c\leq l<r\leq d$ and $x\in[l, r]$,
\begin{equation*}
u(x)\geq
u(l)\frac{F(r)-F(x)}{F(r)-F(l)}+u(r)\frac{F(x)-F(l)}{F(r)-F(l)}.
\end{equation*}

Consider the solution of $(X,Z)$ of the stochastic differential
equation with reflection
\begin{equation}\label{eq:state}
dX_t=\mu(X_t)dt + \sigma(X_t)dW_t-dZ_t, \quad X_{0-}=x\in (c, d),
\end{equation}
where $Z=(Z_{t})_{t \geq 0}$ is a continuous non-decreasing
(except at $t=0$) $\{\mathcal{F}_t\}$-adapted process such that
\begin{equation}\label{eq:Z}
Z_t-Z_0=\int_{(0,t)}1_{\{X_s=b\}}dZ_s,
\end{equation}
for some $b \in (c,d)$. Here, we use the same notation as
\cite{JBS1995}, see equations (4.7) and (4.8). To emphasize the
fact that the initial value of the process $Z$, $Z_0$, depends on
$X_{0-}=x$, below we denote it by $Z_0(x)$. We assume that $x
\rightarrow Z_0(x)$, $x \in (c,d)$, is a measurable function and
$Z_0(c)=0$.

Here, $Z$ is called the \emph{local time} of the process $X$ at
point $b$. When this control is applied to the state process
$\{X_t\}_{t \geq 0}$, for $t>0$, it moves in $(c,b]$ and it is
reflected at $b$ (until the time of absorption). First, we will
find the \emph{best} local time strategy. We will denote the set
of local time strategies by $\CA$. Next, in
Proposition~\ref{prop:more-than-local}, we will show that under
some certain assumptions the local time strategies are optimal in
a larger class of controls, namely non-decreasing,
$\{\mathcal{F}_t\}$-adapted controls. Let $\tau_c \triangleq
\inf\{t \geq 0: X_t=c\}$.
 We consider the following
performance measure associated with $Z \in\CA$
\begin{equation}\label{eq:J}
    J^Z(x)=h \cdot Z_0(x)+\E^{x-Z_0(x)}\left[\int_0^{\tau_c}e^{-\alpha s}
    f(X_s)ds+h \int_{(0,\tau_c)}e^{-\alpha s}dZ_s\right],
\end{equation}
for some given $h \in \R_{+}$. Here, $\P^{x-Z_0(x)}\{\cdot\}$ is a
short-hand notation for the conditional probability measure
$\P\{\cdot|X_0=x-Z_0(x)\}$ and $\E^{x-Z_0(x)}$ is the expectation
with respect to that probability measure. In (\ref{eq:J}), we used
the following notation
\[
h\int_{(0,\tau_c)}e^{-\alpha s}dZ_s \triangleq
h\int_0^{\tau_c}e^{-\alpha s}dZ_s-h Z_0(x).
\]
 The objective is to find the optimal strategy $Z^* \in
\CA$ (if it exists) and the value function:
\begin{equation}\label{problem1}
v(x)\triangleq \sup_{Z\in \CA}J^{Z}(x)=J^{Z^*}(x).
\end{equation}
One could choose $\CA$ to be the family of non-decreasing,
$\{\mathcal{F}_t\}_{ t\geq 0}$-adapted process. We will show in
Proposition~\ref{prop:more-than-local} that it is enough to
consider only the local time strategies under certain assumptions.

\begin{assump}\label{as:local-time}\normalfont
The function $f:(c,d)\rightarrow \mathbb{R}$ is continuous and
satisfies
\begin{equation} \label{eq:f-condition}
\ME\left[\int_0^\infty e^{-\alpha s}|f(X^0_s)|ds\right]<\infty.
\end{equation}
\end{assump}
\subsection{Characterization of the Value Function Corresponding to a Given Reflection Level}

We will first obtain a dynamic programming equation for the
performance measure (\ref{eq:J}). Next, we will apply a
transformation to linearize the difference between the value
associated with a particular control and the value associated with
not applying any control at all. Recall that the region in which
the particular control prescribes no action is commonly referred
to as the continuation region or inaction region of this
particular control.

Let $\tau_b \triangleq\{t \geq 0: X_t \geq b\}$. On denoting
\begin{equation}\label{eq:defn-g}
 g(x)\triangleq \ME \left[\int_0^\infty e^{-\alpha
s}f(X_s^0)ds\right],
\end{equation}
for $x \in [c,b]$, we can write
\begin{align*}
\mathbb{E}^{x-Z_0(x)}&\left[\int_0^{\tau_c}e^{-\alpha
s}f(X_s)ds\right]=\ME\left[\int_0^{\tau_c}e^{-\alpha
s}f(X_s)ds\right]\\ &=\ME\left[\int_0^{\tau_b \wedge \tau_c
}e^{-\alpha s}f(X^0_s)ds+ e^{-\alpha (\tau_b \wedge \tau_c)}
\E^{X_{\tau_b} \wedge \tau_c} \int_{0}^{\tau_c}e^{-\alpha
s}f(X_s)ds\right]\\ &=g(x)-\ME[e^{-\alpha ( \tau_c \wedge \tau_b)
}g(X^0_{\tau_b \wedge \tau_c})]+\ME\left[e^{-\alpha (\tau_b \wedge
\tau_c)}\E^{X_{\tau_b \wedge \tau_c}} \int_{0}^{\tau_c}e^{-\alpha
s}f(X_s)ds\right]\\ &=g(x)-\ME[e^{-\alpha(\tau_b \wedge
\tau_c)}g(X_{(\tau_b \wedge \tau_c)-})]+\ME\left[e^{-\alpha
(\tau_b \wedge \tau_c)}\E^{X_{\tau_b \wedge \tau_c}}
\int_{0}^{\tau_c}e^{-\alpha s}f(X_s)ds\right],
\end{align*}
in which the third line follows from Lemma~\ref{lem:f-app-smp} in
the Appendix. Therefore, for $x \in (c,b]$
\begin{equation}\label{eq:dyn-pro-J}
\begin{split}
J^Z(x)&=\ME\left[e^{-\alpha (\tau_b \wedge
\tau_c)}\left\{-g(X_{(\tau_b \wedge \tau_c) -})+\E^{X_{\tau_b
\wedge \tau_c}} \left[\int_{0}^{\tau_c}e^{-\alpha
s}f(X_s)ds+\int_{(0,\tau_c)} e^{-\alpha s}h \,
dZ_s\right]\right\}\right]+g(x)\\ &=\ME\left[e^{-\alpha(\tau_b
\wedge \tau_c)}\left\{-g(X_{(\tau_b \wedge \tau_c)
-})+J^Z(X_{\tau_b \wedge \tau_c})\right\}\right]+g(x).
\end{split}
\end{equation}
Let us define
\begin{equation}
u^b(x)\triangleq J^Z(x)-g(x), \quad x \in [c,d).
\end{equation}
It is worth noting that $u^b(c)=-g(c)$ since $J^{Z}(c)=0$.

Equation (\ref{eq:dyn-pro-J}) can be written as
\begin{align}\label{eq:u}
u^b(x)&=\ME\left[e^{-\alpha(\tau_b \wedge \tau_c) }\{u^b(X_{\tau_b
\wedge \tau_c})-g(X_{(\tau_b \wedge \tau_c) -})+g(X_{\tau_b \wedge
\tau_c})\}\right]\nonumber\\ &=\ME\left[e^{-\alpha(\tau_b \wedge
\tau_c)}u^b(X_{\tau_b \wedge \tau_c})\right],
\end{align}
for $x \in (c,b]$.   On the other hand, if $x \in [b,d)$, then
\begin{equation}\label{eq:u-lt}
u^b(x)=h\cdot(x-b)-g(x)+g(b)+u^b(b), \quad x \in [b,d).
\end{equation}
Using  (\ref{eq:u}) and (\ref{eq:u-lt}) can be written in a more
compact form as
\begin{eqnarray}\label{eq:singular-twoside}
u^b(x)&=
\begin{cases}
u_0^b(x)\triangleq\ME\left[1_{\{\tau_b<\tau_c\}}e^{-\alpha\tau_b}u^b(b)+1_{\{\tau_b>\tau_c\}}e^{-\alpha\tau_c}u^b(c)\right],
&x\in[c, b],\\K(x, b)+u^b_0(b), &x\in[b, d),
\end{cases}
\end{eqnarray}
in which
\begin{equation}\label{eq:denf-K}
K(x, y)\triangleq h\cdot(x-y)-g(x)+g(y).
\end{equation}
Observe that $u^b(x)$ is continuous at $x=b$.

Using Lemma~\ref{lem:laplace} we can write the function $x
\rightarrow u_0^{b}(x)$, $x \in (c,b]$ as
\begin{equation}\label{eq:ub-c}
u_0^{b}(x)=u^b(b) \frac{\psi(c) \varphi(x)-\psi(x)
\varphi(c)}{\psi(c) \varphi(b)-\psi(b)\varphi(c)}+u^b(c)
\frac{\psi(x)\varphi(b)-\psi(b)
\varphi(x)}{\psi(c)\varphi(b)-\psi(b) \varphi(c)}.
\end{equation}

The function $x \rightarrow u^{b}_0(x)$, $x \in (c,b]$ can be
linearized by using
\begin{equation}
W^b_0(x) \triangleq (u_0^b/\varphi)\circ F^{-1}(x), \quad x \in
[F(c),F(b)],
\end{equation}
and (\ref{eq:ub-c}) becomes
\begin{equation}\label{eq:W0}
W_0^b(x)=W_0^b(F(c))\frac{F(b)-x}{F(b)-F(c)}+W_0^b(F(b))\frac{x-F(c)}{F(b)-F(c)},
\quad x\in [F(c), F(b)].
\end{equation}
We extend the function $x \rightarrow W_0^b(x)$, from $x \in
[F(c),F(b)]$ to $[F(c),F(d))$ by defining
\begin{equation}\label{eq:W}
W^b(x) \triangleq (u^b/\varphi)\circ F^{-1}(x), \quad x
\in[F(c),F(d)).
\end{equation}

We have now established that $W^b(x)$ is a \emph{linear function}
in the transformed \emph{continuation region} (the region of no
action). Note that (\ref{eq:singular-twoside}) and (\ref{eq:W0})
do not completely determine $u^b$: the slope and the intercept of
the line $W^b(x)$ $x \in [F(c),F(b)]$ need to be determined. But
we already know that
\begin{equation}\label{eq:text}
\text{the linear function $W^{b}_0(\cdot)$
 passes through $(F(c),l_c)=\left(F(c),
\frac{-g(c)}{\varphi(c)}\right)$.}
\end{equation}
The slope of this linear function will be determined as a function
of $b$, i.e., $b \rightarrow \beta(b)$, $b \in (c,d)$. Then, we
will give sufficient conditions in
Proposition~\ref{prop:singular-solution} under which the optimal
$b^*$, i.e. $b \in (c,d)$ such that $u^{b^*}(x)=v(x)-g(x)$, can be
determined by the ordinary first order condition, i.e. as the
unique solution of $\partial \beta(b)/\partial b=0$.

\subsection{Characterization of the Optimal Reflection Level}

In this section, we characterize the optimal level at which the
diffusion is to be reflected to maximize a given reward functional
as the unique solution of a non-linear equation. %Moreover, we give
%sufficient conditions for the singular control problem to be
%$\mathbb{C}^{2}$ on its domain.
 \noindent We first transform the function $K(\cdot)$, defined in (\ref{eq:denf-K}), into
\begin{equation}\label{eq:tansform-K}
R(x; b)\triangleq \frac{K(F^{-1}(x), b)}{\varphi(F^{-1}(x))},
\quad x\in [F(b),F(d)).
\end{equation}
From (\ref{eq:W0}) and (\ref{eq:text}) it follows that
\begin{equation}\label{eq:W-with-beta-d}
W^b(x)=\beta (x-F(c)) + l_c, \quad x \in [F(c),F(b)],
\end{equation}
for some $\beta \in \R$, which is to be determined as a function
of $b$. Our task in this section is to identify an appropriate
slope $\beta^*=\beta(b^*)$, so that the function $b \rightarrow
W^{b}(x)$ is maximized at $b^*$ for any $x \in (c,d)$.

\begin{proposition}\label{prop:singular-solution}
Let us define $k: \R \rightarrow \R$ by
\begin{equation}\label{eq:k}
k(x)\triangleq h-g'(x)-l_c\varphi'(x).
\end{equation}
Assume that: (i)\,For any $b \in [c,d)$, $R(y; b)$ defined in
(\ref{eq:tansform-K}) is differentiable with respect to $y$;
(ii)\,For any $b \in [c,d)$, $x \rightarrow R(x;
b)+\frac{\varphi(b)}{\varphi(F^{-1}(x))}W(F(b))$ is increasing and
concave on $x\in (F(j), F(d))$ for some point $j \in (c, d)$ and
it approaches infinity as $x \rightarrow d$; (iii)\,There exists a
unique solution $b^* \in (c,d)$ to the equation
\begin{equation}\label{eq:optimal-b}
k' (b)\psi'(b) -k(b) \psi''(b)+F(c)[k (b)\varphi''(b)-k' (b)
\varphi'(b)]=0
%(F''(b)\varphi(b)+2F'(b)\varphi'(b)+F(b)\varphi''(b))=(g''(b)+d\varphi''(b))(F'(b)\varphi(b)+F(b)\varphi'(b))
\end{equation}
such that $b^*$ satisfies
\begin{equation}\label{eq:sufficiency}
k''(b^*)\psi'(b^*)-k(b^*)\psi'''(b^*)+F(c)(-k''(b)\varphi'(b)+k(b)
\varphi''(b))<0.
\end{equation}
Then the solution of (\ref{problem1}) is given by
$v(x)=u^{b^*}(x)+g(x)$, in which $u^{b^*}$ is given by the
equation (\ref{eq:singular-twoside}) if we replace $b$ by $b^*$
and choose the slope of (\ref{eq:W0}) to be $\beta^*$, which is
given by
\begin{equation}\label{eq:optimal-beta}
\beta^* \triangleq \frac{k(b^*)}{\psi'(b^*)-F(c)\varphi'(b^*)}.
\end{equation}
Recall that $W^{b^*}(x)$ is linear for $x \in [F(c),F(b^*))$.
\end{proposition}

\begin{proof}
We will first determine the slope of the line in
(\ref{eq:W-with-beta-d}), as a function of $b$, i.e., $b
\rightarrow \beta(b)$, $b \in [c,d)$. This will be established by
showing that $W^b$ defined in (\ref{eq:W}) is continuously
differentiable at $b$. To this end we will first consider the
threshold strategy that is characterized by the pair $(b,a) \in
[c,d)^2$: $Z$ is said to be a threshold strategy corresponding to
$(b,a)$ if, whenever the process $X$ in (\ref{eq:state}) hits
level $a$ or is above $a$, then it jumps to (the jump is forced by
$Z$) level $b \leq a$. (Although, the letter $Z$ was used to
denote only local time strategies before we would like to use it
to denote the threshold strategies to be able to refer
(\ref{eq:state}) when we are describing threshold strategies. This
prevents introducing unnecessary equations.) Consider the reward
in (\ref{eq:J}) corresponding to the particular threshold strategy
$Z$ and denote it by $u^{b,a}$. Note that $u^{b,b}=u^b$. The
control represented by the pair $(b,a)$ is of impulse control
type. For $b \in [c,d)$, let us find $a(b)$ such that $\sup_{a \in
[c,d)}u^{b,a}=u^{b,a(b)}$. Using the results in Section 2.2 of
\cite{BE2006} and the assumption (i) and (ii) of the proposition
we conclude that, for any $b$, there exists a unique $a(b)$ such
that the function $x \rightarrow W^{b, a(b)}(x)$, defined by
(\ref{eq:W}) when $u$ is replaced by $u^{b, a(b)}$, is
continuously differentiable at $a(b)$. This characterization of
the function $b \rightarrow a(b)$, $b \in [c,d)$ will be used to
show $a(b)=b$, $b \in [c,d)$ and to calculate the slope $b
\rightarrow \beta(b)$, $b \in [c,d)$.

Let us define
\begin{equation}\label{eq:W-tilde-b}
\tilde{W}^b(x) \triangleq \begin{cases} \beta(b) (x-F(c)) + l_c &
x \in [F(c),F(b)]
\\ H(x,b)\triangleq R(x; b)+\frac{\varphi(b)}{\varphi(F^{-1}(x))}\left(\beta(b)(F(b)-F(c))+l_c\right), & x
\in[F(b),F(d)).
\end{cases}
\end{equation}
in which
\begin{equation}\label{eq:betab}
\beta(b)=\frac{k(b)}{\psi'(b)-F(c) \varphi'(b)}.
\end{equation}
The right-hand derivative of the function $\tilde{W}^b$ satisfies
\begin{equation}\label{eq:rem2}
\begin{split}
&(\tilde{W}^{b})'(F(b)+)=-(\beta(b)
(F(b)-F(c))+l_c)\frac{1}{\varphi(b)}\frac{\varphi'(b)}{F'(b)}+\frac{\partial}{\partial
x}R(x;b)\bigg|_{x=F(b)}
\\&=-(\beta(b)
(F(b)-F(c))+l_c)\frac{1}{\varphi(b)}\frac{\varphi'(b)}{F'(b)}+\frac{\varphi(x)(-g'(x)+h)-(h
\cdot
(x-b)-g(x)+g(b))\varphi'(x)}{\varphi(x)^2}\frac{1}{F'(x)}\bigg|_{x=b}
\\& =\beta(b),
\end{split}
\end{equation}
where we used (\ref{eq:denf-K}) and (\ref{eq:tansform-K}) to
derive the first equality, and (\ref{eq:betab}) to derive the
third inequality. The function $x \rightarrow \tilde{W}^{b}(x)$,
$x \in [F(c),F(d))$ is $C^1$ at $x=F(b)$, since the left-hand
derivative is also $\beta(b)$. This implies that
$\tilde{W}^b=W^{b, a(b)}$ or $(u^{b,a(b)}\varphi) \circ
F^{-1}=\tilde{W}^{b}$ and that $a(b)=b$. As a result we see that
$W^b$ satisfies smooth fit condition at $b$ and the slope in
(\ref{eq:W-with-beta-d}) is given by (\ref{eq:betab}). Before we
continue with the proof the reader should note that
\begin{align*}
\beta(b)&=\lim_{a\downarrow b}\frac{R(F(a); b))+l_c(
\frac{\varphi(b)}{\varphi(a)}-1)}{F(a)
-\frac{\varphi(b)}{\varphi(a)}F(b)+F(c)\left(-1+\frac{\varphi(b)}{\varphi(a)}\right)}=\lim_{a\downarrow
b}\frac{h-g'(a)-l_c
\varphi'(a)}{F'(a)\varphi(a)+F(a)\varphi'(a)-F(c)\varphi'(a)}\\
&=\frac{h-g'(b)-l_c
\varphi'(b)}{F'(b)\varphi(b)+F(b)\varphi'(b)-F(c)\varphi'(b)}=\frac{k(b)}{\psi'(b)-F(c)\varphi'(b)},
\end{align*}
where the second equality follows from an application of
L'Hospital's rule. In contrast with equation (2.25) in
\cite{BE2006}, this implies that the first order smooth fit of the
singular control at $b$ can be derived by a limiting argument from
the continuous fit of a family of impulse control problems at $b$.
Here, the first order smooth fit holds at any $b \in (c,d)$, not
only at $b^*$, which we will soon discover to be the optimal
reflection barrier.

Equations (\ref{eq:optimal-b}) and (\ref{eq:sufficiency}) imply
that $\beta'(b^*)=0$ and $\beta''(b^*)<0$. Therefore, $b^*$ is a
local maximum of the function $b \rightarrow \beta(b)$. On the
other hand, since we assumed that the uniqueness of the solution
to (\ref{eq:optimal-b}), $b^*$ is the unique local extremum of the
function  $b \rightarrow \beta(b)$. Let us argue that $b^*$ is the
global maximum of this function: Assume there exists a point $m
\neq b^*$ where the maximum of the function $b \rightarrow
\beta(b)$ is attained. Then there would exist local minimum $n \in
(m,b^*)$ of the function $b \rightarrow \beta(b)$ which
contradicts the fact that $b^*$ is the unique local extremum of
this function. Note that $b \rightarrow \beta(b)=k(b)/\psi'(b)$
may not be concave.

%Observe that $b^*$ defined in the statement of proposition is the
%unique local extremum of the function $\beta(b)=k(b)/\psi'(b)$,
%and moreover it is a local maximum. This implies that $b^*$ is the
%global maximum of $\beta(b)$. If the maximum were to be attained
%at the end points of the interval $(c,d)$, then there would be a
%local minimum since we already know that $b^*$ is a local maximum.
%But this contradicts the fact that $b^*$ is the only local
%extremum of $\beta(b)$.

%Now, after some algebraic manipulations using $\beta'(b^*)=0$ and
%$\beta''(b^*)<0$,
Recall the definition of the function $H$ from
(\ref{eq:W-tilde-b}). Using (\ref{eq:denf-K}) and
(\ref{eq:tansform-K}) we can calculate the derivative of $H$ with
respect to $b$ as
\begin{equation}\label{eq:der-H}
\begin{split}
\frac{\partial}{\partial
b}H(x;b)\bigg|_{b=b^*}&=\frac{-h+g'(b)+\varphi'(b)(\beta(b)
(F(b)-F(c))+l_c
)+(\beta'(b)(F(b)-F(c))+\beta(b)F'(b))\varphi(b)}{\varphi(F^{-1}(x))}\bigg|_{b=b^*}
\\&=\frac{-h+g'(b)+\beta(b)\psi'(b)+l_c \varphi'(b)+\beta'(b)F(b)\varphi(b)-F(c)(\beta(b)\varphi'(b)+\beta'(b)\varphi(b))}{\varphi(F^{-1}(x))}\bigg|_{b=b^*}
\\&=\frac{\beta'(b)(F(b)-F(c))\varphi(b)}{\varphi(F^{-1}(x))}\bigg|_{b=b^*}=0.
\end{split}
\end{equation}
Here, the second equality follows from the definition of $F$ in
(\ref{eq:F}), and the third equality follows from the definition
of $\beta(b)$ in (\ref{eq:betab}) and (\ref{eq:k}). Note that,
\begin{equation}\label{eq:iff}
\frac{\partial}{\partial b}H(x;b)=0 \quad \text{if and only if
$b=b^*$}.
\end{equation}
On the other hand,
\begin{equation}\label{eq:sec-der-H}
\frac{\partial^2}{\partial
b^2}H(x;b)\bigg|_{b=b^*}=\frac{\beta''(b)(F(b)-F(c))\varphi(b)+\beta'(b)F'(b)\varphi(b)+\beta'(b)(F(b)-F(c))\varphi'(b)}
{\varphi(F^{-1}(x))}\bigg|_{b=b^*}<0,
\end{equation}
since $\beta'(b^*)=0$, $\beta''(b^*)<0$ and $F$ is increasing.
Now, (\ref{eq:W-tilde-b}), (\ref{eq:der-H}), (\ref{eq:iff}) and
(\ref{eq:sec-der-H}), together with the fact that $b \rightarrow
\beta(b)$ is maximized at $b^*$ imply that
\begin{equation}\label{eq:sup-W}
W^{b*}(x)=\sup_{b \in (c,d)}W^{b}(x), \quad x \in [F(c),F(d)).
\end{equation}
The proof of our assertion follows since it is immediate from
(\ref{eq:sup-W}) that $u^{b^*}(x)=\sup_{b \in [c,d)}u^b(x)$, for
all $x \in [c,d)$.
\end{proof}

We can extend our results to solve constrained optimization
problems. A dividend payout problem with solvency constraints was
recently solved by Paulsen \cite{paulsen}. In this problem, the
firm is allowed to pay dividends only if the cash flow process $X$
is greater than certain (pre-determined) value $\tilde{b}$. Here,
we provide a simpler proof to this result, using the
characterization we provided for the optimal reflection barrier in
Proposition~\ref{prop:singular-solution}.
\begin{corollary} Assume that the assumptions of
Proposition~\ref{prop:singular-solution} hold. Let $\tilde{\CA}$
be the set of $Z \in \CA$ such that
\begin{equation}\label{eq:Z-p}
Z_t-Z_0=\int_{(0,t)}1_{\{X_s=b\}}dZ_s, \quad b\leq \tilde{b} \in
(c,d),
\end{equation}
for a fixed $\tilde{b}$ and
 define
\begin{equation}
\hat{b} \triangleq \begin{cases} b^* & \text{if} \quad b^* \leq
\tilde{b},
\\ \tilde{b} & \text{if} \quad b^*> \tilde{b}.
\end{cases}
\end{equation}
Let us also define $\hat{Z}$ by replacing $b$ with $\hat{b}$ in
(\ref{eq:Z-p}). Then
\begin{equation}
v(x)\triangleq \sup_{Z\in \tilde{\CA}}J^{Z}(x)=J^{\hat{Z}}(x).
\end{equation}
\end{corollary}

\begin{proof}
The proof follows since under the assumptions of
Proposition~\ref{prop:singular-solution} the functions $b
\rightarrow \beta(b)$, $b \in [c,d)$, $b \rightarrow H(x;b)$, $x
\in [c,d)$ have a unique maximum at $b^*$. (See the proof of
Proposition~\ref{prop:singular-solution}).
\end{proof}

\subsection{The Optimality of Local Time
Strategies in the Class of Monotone Increasing
Controls}\label{sec:optimality}

Let us write the value function $v(x)$, explicitly and make some
observations on it.
\begin{eqnarray}\label{eq:solution}
v(x)&=
\begin{cases}
v_0(x)\triangleq \varphi(x)(\beta^*(F(x)-F(c))+l_c)+g(x), &x\in[c,
b^*],\\ h \cdot (x-b^*)+v_0(b^*), &x\in[b^*, d)
\end{cases}
\end{eqnarray}
where the second equation is obtained by
\begin{align*}
K(x, b^*)+ u^{b^*}_0(b^*)+g(x)&=h \cdot
(x-b^*)-g(x)+g(b^*)+\varphi(b^*)(\beta^*(F(b^*)-F(c))+l_c)+g(x)\\
&=h\cdot(x-b^*)+v_0(b^*).
\end{align*}
It is worth noting that $v(c)=0$. \noindent
\begin{remark}\normalfont \label{rem:verification-remark}
\begin{enumerate}
\item [(a)]
The first and the second derivative of $v(x)$ on $(c,b^*)$ are
\begin{align*}
v'(x)&=\beta^*\psi'(x)+(l_c-\beta^* F(c))\varphi'(x)+g'(x)\quad
\text{and}\quad v''(x)=\beta^*\psi''(x)+(l_c-\beta^*
F(c))\varphi''(x)+g''(x).
%\frac{h-d\varphi(b^*)}{F'(b^*)\varphi'(b^*)+F(b^*)\varphi'(b^*)}(F'(b^*)
%\varphi(b^*)+F(b^*)\varphi'(b^*))+d\varphi(b^*)+g'(b^*)\\
%&=h+g'(b^*).
\end{align*}
Evaluating these expressions at $b^*$ we obtain
\begin{equation} \label{eq:first-second-smooth-fit}
\begin{split}
v'(x)\bigg|_{x=b^*}&=\frac{k(b^*)}{\psi'(b^*)-F(c)\varphi'(b^*)}(\psi'(x)-
F(c) \varphi'(x) )+h-k(x)\bigg|_{x=b^*}=h,
\\ v''(x)\bigg|_{x=b^*}&=\frac{k(b^*)}{\psi'(b^*)-F(c)\varphi'(b^*)}(\psi''(x)-F(c) \varphi''(x))-k'(x) \bigg|_{x=b^*}=0.
\end{split}
\end{equation}
We used  (\ref{eq:k}) and (\ref{eq:optimal-beta}) to obtain the
first expression and (\ref{eq:optimal-b}) to obtain the second
expression. Note that these smooth fit conditions are the two
boundary conditions that are frequently imposed to solve the
singular control problems in an ordinary Hamilton-Jacobi-Bellman
(HJB) approach. In that approach, after the solution is
constructed, the assumptions are verified using a verification
lemma. However, the smooth fit conditions need not necessarily
hold and the HJB approach is unable to tell the sufficient
conditions for the smooth fit to hold. Using our alternative
methodology, in Proposition~\ref{prop:singular-solution}, we are
able to list some sufficient conditions for the value function to
be $\mathbf{C}^{2}$ on all of its domain.

Furthermore,
\begin{equation}\label{eq:infinitesimal-equality}
(\mathcal{A}-\alpha)v(x)=(\mathcal{A}-\alpha) g(x)= -f(x) \quad
\text{for} \quad x \in (c,b^*].
\end{equation}
\item[(b)] Under assumption (iii) of
Propostion~\ref{prop:singular-solution} the function $b
\rightarrow k(b)/(\psi'(b)-F(c)\varphi'(b))$ is maximized at
$b^*$. Therefore, using the first equation in
(\ref{eq:first-second-smooth-fit}) it can be checked that $v'(x)
\geq h$, $x \in (c,d)$.

\item [(c)] Since $b^*$ satisfies (\ref{eq:optimal-b}), then we have that
\begin{equation}
\frac{k(b^*)}{\psi'(b^*)-F(c)\varphi'(b^*)}=\frac{k'(b^*)}{\psi''(b^*)-F(c)\varphi''(b^*)}.
\end{equation}
Now using, the second equation in
(\ref{eq:first-second-smooth-fit}) we have that
\begin{equation}
v''(x)=\frac{k'(b^*)}{\psi''(b^*)-F(c)\varphi''(b^*)}(\psi''(x)-F(c)\varphi''(x))-k'(x),
\quad x \in  (c,b^*]
\end{equation}
If either
\begin{equation} \label{eq:sec-der-1}
\varphi''(x)-F(c)\varphi''(x)>0, \quad x \in (c,b^*]\quad
\text{and} \quad x \rightarrow \frac{k'(x)}{\psi''(x)-F(c)
\varphi''(x)} \quad \text{is a decreasing function on $(c,b^*]$},
\end{equation}
or
\begin{equation}\label{eq:sec-der-2}
\varphi''(x)-F(c)\varphi''(x)<0, \quad x \in (c,b^*]\quad
\text{and} \quad x \rightarrow \frac{k'(x)}{\psi''(x)-F(c)
\varphi''(x)} \quad \text{is an increasing function on $(c,b^*]$},
\end{equation}
then $v''(x) \leq 0$, $x \in (c,b^*].$ Note from
(\ref{eq:solution}) that $v''(x)=0$ on $x\in [b^*, d)$.

\item [(d)]For $x \in [b^*,d)$
\begin{equation}\label{eq:infinitesimal-inequality}
\begin{split}
(\mathcal{A}-\alpha)v(x)&= (\mathcal{A}-\alpha)
(h\cdot(x-b^*)+v_0(b^*))=\mu(x) h -\alpha h\cdot(x-b^*)-\alpha
v_0(b^*),
\\ &\leq \mu(x) h-\alpha v_0(b^*) \leq \mu(b^*)  h-\alpha v_{0}(b^*)
=\lim_{x \downarrow b^*} (\mathcal{A}-\alpha)v(x)
\\ &=\lim_{x \uparrow b^*} (\mathcal{A}-\alpha)v(x)=-\lim_{x
\uparrow b^*} f(x)=-f(b^*) \leq -f(x)
\end{split}
\end{equation}
if we assume that the maximums of the functions $x \rightarrow
\mu(x)$ and $x \rightarrow f(x)$ on the interval $[b^*,d)$ are
attained at $b^*$ (for e.g. if both $x \rightarrow f(x)$ and $x
\rightarrow \mu(x)$ are non-increasing on $[b^*,d)$), and that
these functions are both continuous at $b^*$.
 Note that the identity $\lim_{x \downarrow b^*}
(\mathcal{A}-\alpha)v(x)=\lim_{x \uparrow b^*}
(\mathcal{A}-\alpha)v(x)$ is due to
(\ref{eq:first-second-smooth-fit}).

\end{enumerate}
\end{remark}

The following proposition gives sufficient conditions under which
the  local time strategies are optimal in the class of all
increasing strategies.
\begin{proposition}\label{prop:more-than-local}
Assume that the assumptions of
Proposition~\ref{prop:singular-solution} hold. Consider the
process
\[ dX_t=\mu(X_t)dt + \sigma(X_t)dW_t-d\xi_t
\]
in which $\xi_t$ is an $\{\mathcal{F}_t\}$-adapted, non-decreasing
and right-continuous process (except possibly at zero) such that
$\E^x\left[\int_{0}^{\infty} e^{-\alpha s} d \xi_s\right]<\infty$.
We denote the family of such controls by $\mathcal{C}$. Let us
assume that $x \rightarrow \sigma(x)$, $x \in (c,d)$ is a bounded
function, the maximums of the functions $x \rightarrow \mu(x)$, $x
\rightarrow f(x)$ on $x \in [b^*,d)$ is attained st $b^*$, and
that both $\mu(\cdot)$ and $f(\cdot)$ are continuous at $b^*$. We
further assume that either (\ref{eq:sec-der-1}) or
(\ref{eq:sec-der-2}) holds. Then $ x \rightarrow v(x)$, $x \in
(c,d)$ defined in (\ref{eq:solution}) satisfies $v(x) \geq
J^{\xi}(x)$, $x \in (c,d)$, for any $\xi \in \mathcal{C}$, in
which
\begin{equation}
    J^{\xi}(x)=h \xi_0(x)+\E^{x-\xi_0(x)}
    \left[\int_0^{\infty}e^{-\alpha s}f(X_s)ds+\int_{(0,\infty)}e^{-\alpha
    s}h d\xi_s\right].
\end{equation}
\end{proposition}
\begin{proof}
We first apply It\^{o}'s formula to $e^{-\alpha t}v(X_t)$ and get
\begin{equation}\label{eq:ito}
\begin{split}
e^{-\alpha t}v(X_t)&=v(x)+\int_{0}^{t} e^{-\alpha
s}(\mathcal{A}-\alpha)v(X_s)ds-\int_{0}^{t}e^{-\alpha
s}v'(X_{s-})d\xi_s +\int_0^{t}e^{-\alpha s}\sigma(X_s)v'(X_s)dW_s
\\&+\sum_{0<s \leq t}e^{-\alpha s}(v(X_s)-v(X_{s-})+v'(X_{s-})\Delta
X_s),
\end{split}
\end{equation}
in which $\Delta X_s=X_{s-}-X_s$, $s \geq 0$. From equation
(\ref{eq:ito}) and Remark~\ref{rem:verification-remark} (b), (c)
and (e), it follows that
\begin{equation}\label{eq:solve-for-v}
\begin{split}
 v(x)&=\int_{0}^{t}h e^{-\alpha s} d \xi_s-\int_{0}^{t}e^{-\alpha
 s}(h-v'(X_{s-}))d\xi_s -\int_0^{t} e^{-\alpha s}
 (\mathcal{A}-\alpha) v(X_s)ds
 \\ &-\int_0^{t}\sigma(X_s)e^{-\alpha s}v'(X_s)dW_s - \sum_{0<s \leq t}e^{-\alpha
s}(v(X_s)-v(X_{s-})+v'(X_{s-})\Delta X_s)+e^{-\alpha t}v(X_t)
\\ & \geq \int_{0}^{t}h e^{-\alpha s} d \xi_s +\int_0^{t} e^{-\alpha s} f(X_s)ds
-\int_0^{t}\sigma(X_s)e^{-\alpha s}v'(X_s)dW_s
\\ & - \sum_{0<s \leq t}e^{-\alpha
s}(v(X_s)-v(X_{s-})+v'(X_{s-})\Delta X_s)+e^{-\alpha t}v(X_t)
\\ & \geq \int_{0}^{t}h e^{-\alpha s} d \xi_s +\int_0^{t} e^{-\alpha s} f(X_s)ds
-\int_0^{t}\sigma(X_s)e^{-\alpha s}v'(X_s)dW_s \\ &- \sum_{0<s
\leq t}e^{-\alpha s}(v(X_s)-v(X_{s-})+v'(X_{s-})\Delta X_s).
 \end{split}
\end{equation}
The last line follows because $v$ is positive: $v(c)=0$ and $v'(x)
\geq h \geq 0$, $x \in (c,d)$, by
Remark~\ref{rem:verification-remark} (b).

We have that
\begin{equation}\label{eq:martingale}
\E\left[\int_0^{t}e^{-\alpha s}\sigma(X_s)v'(X_{s})dW_s \right]=0,
\end{equation}
since $v'(x)$ and $\sigma(x)$ are bounded. On the other hand,
since $v''(x) \leq  0$ (see Remark~\ref{rem:verification-remark}
(d)) for any $x>y$
\[
v(x)-v(y)-v'(x)(x-y)=\int_{y}^{x}(v'(u)-v'(x))du \geq 0,
\]
which implies that
\begin{equation}\label{eq:jump-term}
- \sum_{0<s \leq t}e^{-\alpha s}(v(X_s)-v(X_{s-})+v'(X_{s-})\Delta
X_s) \geq 0.
\end{equation}

 Now (\ref{eq:solve-for-v}), (\ref{eq:martingale}),
(\ref{eq:jump-term}) imply that
\[
v(x) \geq \E^x\left[\int_{0}^{t}h \cdot e^{-\alpha s} d \xi_s
+\int_0^{t} e^{-\alpha s} f(X_s)ds\right],
\]
for all $t>0$, which implies that $v(x) \geq J^{\xi}(x)$ for all
$x \in (c,d)$ after taking a limit as $t\rightarrow \infty$. The
exchange of limit and integration is possible due to
Assumption~\ref{as:local-time} and the definition of $\mathcal{C}$
as a result of an application of bounded convergence theorem.
\end{proof}

\begin{remark}\label{rem:check}
\normalfont We give a useful hint which will be helpful in
checking whether
\\ $x \rightarrow R(x;
b)+\frac{\varphi(b)}{\varphi(F^{-1}(x))}W(F(b))$, $x \in (F(c),
F(d))$ satisfies assumption (ii) of
Proposition~\ref{prop:singular-solution}. Let us denote
\begin{equation}\label{eq:m(x)}
m(x)= \frac{1}{F^{'}(x)}\left(\frac{K}{\varphi}\right)^{'}(x),
\end{equation}
then $R^{'}(y;b)=m(x)$ and $R^{''}(y;b)=m^{'}(x)/F^{'}(x)$, in
which $y\triangleq F(x)$.
 If $x \rightarrow K(x,b)$ is twice-differentiable at $x \in
(c,d)$, then
\begin{equation}\label{eq:devH}
\quad R^{''}(y;b) [(\mathcal{A}-\alpha)K(x,b)]\geq 0.
\end{equation}
The inequality is strict if $R^{''}(y;b)\neq 0$.
\end{remark}

\subsection{Examples of Dividend Payment Problems}

\begin{example} \normalfont
\textbf{Dividend payout with a Brownian motion with drift
(Jeanblanc and Shiryaev \cite{JBS1995}, Case C):} Let us assume
that the capital of a company is modeled a Brownian motion with a
drift and the managers of the company would like to maximize the
amount of dividends payed out. We assume that the company is
ruined when the capital becomes $0$ (i.e. 0 an absorbing
boundary). The right boundary $+\infty$ is natural.  The
uncontrolled process $X^0$ is a Brownian motion with drift
\begin{equation*}
dX^0_t=\mu dt +\sigma dW_t.
\end{equation*}  The value function is defined as
\begin{equation}\label{eq:JS}
v(x)\triangleq\sup_{\CA}\left\{Z_0(x)+\E^{x-Z_0(x)}\left[\int_{(0,
\tau_0)}e^{-\alpha t}dZ_t\right]\right\},
\end{equation}
where $\tau_0=\inf\{t\le 0; X_t=0\}$.

In this problem, $f(x)\equiv 0$ and $h\equiv 1$ and $K(x, y)=x-y$.
 As in \cite{JBS1995} we take $\sigma=\sqrt{2}$. By solving the
equation $(\mathcal{A}-\alpha)v(x)=0$, in which $\mathcal{A}$ is
the infinitesimal generator of the uncontrolled process $X^0$, we
find $\psi(x)=e^{(-\frac{\mu}{2}+\Delta)x}$ and
$\varphi(x)=e^{(-\frac{\mu}{2}-\Delta)x}$ where
$\Delta=\sqrt{(\frac{\mu}{2})^2+\alpha}$.  Hence $F(x)=e^{2\Delta
x}$ and $F^{-1}(x)=\frac{\log x}{2\Delta}$.  Note that $F(0)=1$,
$l_c=0$ and $k(x)=1$.

\noindent \textbf{Verification of the Conditions in Proposition
\ref{prop:singular-solution}.}

\noindent (i) For a given $b>0$, we have that
\begin{equation*}
R(y; b)=(K/\varphi)(F^{-1}(y))=\frac{\log
y}{2\Delta}\left(y^{\frac{\frac{1}{2}\mu+\Delta}{2\Delta}}-b\right)
\end{equation*}
on $y>F(0)=1$.  This function is differentiable with respect to
$y$.\\ (ii) $R(\cdot; b)$ is increasing on $y\in [F(0), \infty)$
by (\ref{eq:m(x)}) and $\lim_{y \rightarrow \infty}
R(y;b)=\infty$. We also have that
\begin{equation}\label{eq:inv-varphi}
\frac{1}{\varphi(F^{-1}(y))}=\frac{\log y}{2\Delta}
\end{equation}
is increasing on $[F(0),\infty)$ to $\infty$.

On the other hand, $(\mathcal{A}-\alpha)K(x, b)=p(x)$ for every
$x>0$, in which $p(x)\triangleq\mu-\alpha(x-b)$. This linear
function $p(x)$ has only one positive root at say, $k$. Then by
(\ref{eq:devH}), $R(y ; b)$ is convex on $y\in [F(0), F(k))$ and
concave on $y\in (F(k), \infty)$.  Observe from
(\ref{eq:inv-varphi}) that $1/\varphi(F^-1(y))$ is concave.

\noindent (iii) From (\ref{eq:optimal-b}), Since $k'(x)=0$ for
$x\in [0, \infty)$, Equations (\ref{eq:optimal-b}) and
(\ref{eq:sufficiency}) become
\begin{equation}\label{eq:JS-bstar}
\psi''(b)=F(0)\varphi''(b),
\end{equation}
\begin{equation}\label{eq:JS-suf}
\psi'''(b^*)-\varphi'''(b^*)F(0)>0.
\end{equation}

Since $\psi''(\cdot)$ is increasing and $\varphi''(\cdot)$ is
decreasing on $[0,\infty)$, equation (\ref{eq:JS-suf}) holds for
all $x\in [0, \infty)$. Moreover, since $\psi''(0)<\varphi''(0)$
and $0=\lim_{x \rightarrow \infty}\varphi''(x)<\lim_{x \rightarrow
\infty} \psi''(x)=\infty$, there exists a unique solution to
(\ref{eq:JS-bstar}).

\noindent \textbf{Verification of the Conditions in Proposition
\ref{prop:more-than-local}:  } The only non-trivial condition to
check is whether $v''(x)\le 0$ for $x\in (0, \infty)$. It can be
shown that $\psi''(x)-F(0)\varphi''(x)< 0$ on $x \in (0, b^*)$, by
the same argument that we used to prove the uniqueness of the root
of (\ref{eq:JS-bstar}) and the concavity of $v$ follows from
Remark~\ref{rem:verification-remark}-c. Hence we conclude that the
local time strategy at $b^*$ is optimal among all the admissible
strategies.

Now, $v_0(\cdot)$, defined in (\ref{eq:solution}) can be computed
as
\begin{align*}
v(x)&=\varphi(x)W^{b^*}(F(x))=\beta^*(F(x)-1)\varphi(x)=\beta^*(e^{2x\Delta}-1)e^{-(\mu/2+\Delta)x}\\
&=\beta^* e^{-\mu x/2}(e^{\Delta x}-e^{-\Delta x})=2\beta^*
e^{-\mu x/2}\sinh(x\Delta).
\end{align*}
  The solution to this problem is then
\begin{eqnarray*}
v(x)&=& \begin{cases}
                    v_0(x), &0\leq x \leq b^*, \\
                    v_0(b^*)+x-b^*, &b^*\leq x,
        \end{cases}
\end{eqnarray*}
which coincides with the solution that is computed by Jeanblanc
and Shiryaev \cite{JBS1995} by using the ordinary HJB approach,
which is specific to the modeling assumptions.
 Figure~\ref{fig:JS} shows the value function after applying the
 transformation (\ref{eq:W}), the slope function $b \rightarrow
 \beta(b)$, $b \in (c,d)$, the value function and its derivative
when the parameters are $(\mu, \alpha)=(0.15, 0.2)$. The optimal
reflection point
 is $b^*=0.736246$ and $\beta(b^*)=1.16523$.
\begin{figure}[h]
\begin{center}
\begin{minipage}{0.45\textwidth}
\centering \includegraphics[scale=0.70]{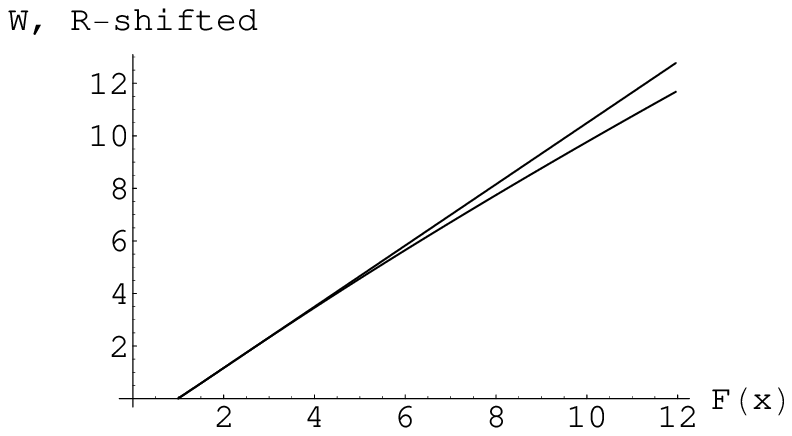} \\ (a)
\end{minipage}
\begin{minipage}{0.45\textwidth}
\centering \includegraphics[scale=0.70]{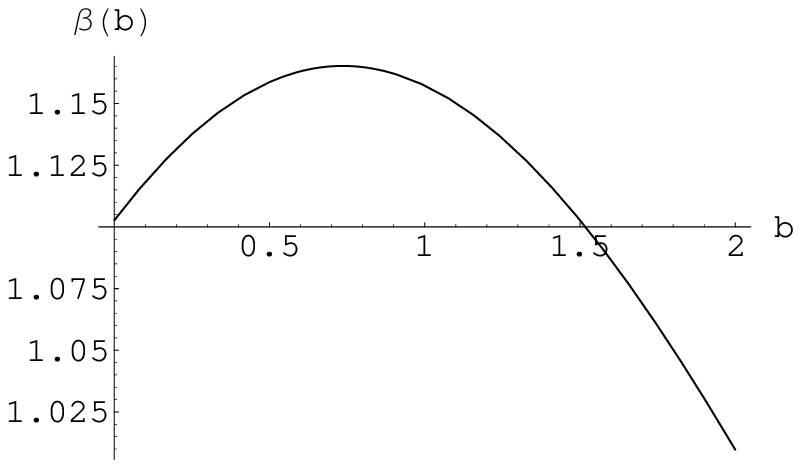} \\ (b)
\end{minipage}
\begin{minipage}{0.45\textwidth}
\centering \includegraphics[scale=0.60]{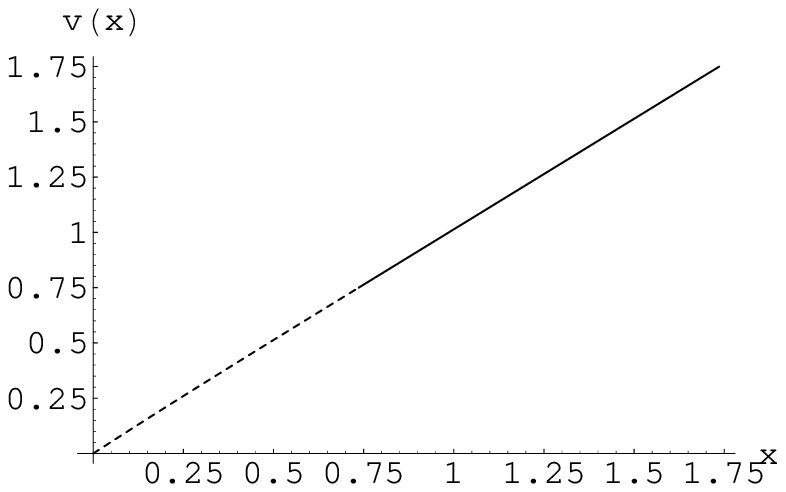} \\ (c)
\end{minipage}
\begin{minipage}{0.45\textwidth}
\centering \includegraphics[scale=0.60]{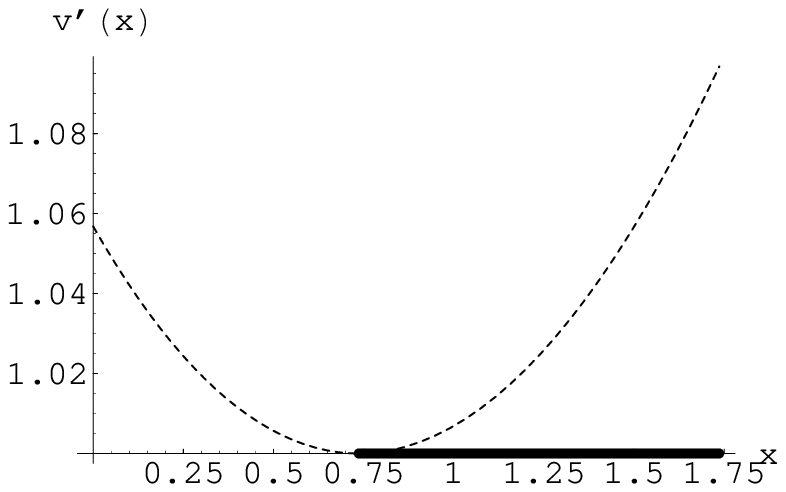} \\ (d)
\end{minipage}
%\begin{minipage}{0.45\textwidth}
%\centering \includegraphics[scale=0.75]{singular4.eps} \\
%(d)
%\end{minipage}
\caption{\small The analysis of the singular control problem of
Jeanblanc and Shiryaev \cite{JBS1995} with parameters
 $(\mu, \alpha)=(0.15, 0.2)$: (a) The function $x \rightarrow R(F(x), b^*)+\varphi(b^*)/\varphi(x)$, $x \in (c,d)$
 and its linear majorant. (b) The graph of $b \rightarrow \beta(b)$ (see (\ref{eq:betab})).
 It attains its maximum at $b^*$. (c) The value function $v(x)$ with $b^*=0.736246$ and $\beta^*=1.16523$.
 (d) The derivative $v'(x)$:  It is  $v_0'(x)$ on $0\leq x \leq b^*$ and $v'(x)=1$ on $b^*\leq x$ (lower line).
 The derivative  $v_0'(x)$ on $b^*\leq x$ is also shown  to illustrate that $v_0''(b^*)=0$.}
\label{fig:JS}
\end{center}
\end{figure}
\end{example}
\begin{example}\normalfont
\textbf{Dividend payout with a square root process:  } We solve
the problem defined in (\ref{eq:JS}) when the cash flow of the
company is modeled by the following square root process:
\begin{equation}\label{eq:CIR}
dX^0_t=(1-2 \rho X_t^0)dt+ 2 \sqrt{X^0_t} dW_t, \quad X^0_0=x>0.
\end{equation}
The solutions of $(\A-\alpha)v(x)=0$ are
\begin{equation} \label{eq:SR-fundamental}
\psi(x)=x^{-1/4}\exp\left(\frac{\rho
x}{2}\right)M_{-\frac{\alpha}{2
\rho}+\frac{1}{4},-\frac{1}{4}}(\rho x), \quad
\varphi(x)=x^{-1/4}\exp\left(\frac{\rho
x}{2}\right)W_{-\frac{\alpha}{2
\rho}+\frac{1}{4},-\frac{1}{4}}(\rho x),
\end{equation}
in which $W_{-\frac{\alpha}{2 \rho}+\frac{1}{4},-\frac{1}{4}}$ and
$M_{-\frac{\alpha}{2 \rho}+\frac{1}{4},-\frac{1}{4}}$ are
Whittaker functions. (See Appendix 1.26 of Borodin and Salminen
\cite{salminen} as well as Chapter 2.1.11.) The Whittaker
functions are defined as
\begin{equation}\label{eq:Whittaker}
\begin{split}
W_{-\frac{\alpha}{2
\rho}+\frac{1}{4},-\frac{1}{4}}\left(\frac{x^2}{2}\right)&=2^{\frac{\alpha}{2
\rho}-\frac{1}{4}}\sqrt{x}D_{-\alpha/\rho}(x), \quad x \geq 0,
\\M_{-\frac{\alpha}{2
\rho}+\frac{1}{4},-\frac{1}{4}}\left(\frac{x^2}{2}\right)&=\frac{\Gamma((1+\alpha/\rho)/2)}{2
\sqrt{\pi}}\sqrt{x}(D_{-\alpha/\rho}(-x)-D_{-\alpha/\rho}(x)),
\quad x \geq 0,
\end{split}
\end{equation}
in which $\Gamma$ stands for the Gamma function
$
\Gamma(x)=\int_{0}^{\infty}u^{x-1}e^{-u}du $ and $D_{\nu}(\cdot)$
is the parabolic cylinder function, which is defined as
\[
D_{\nu}(x)\triangleq
2^{\nu/2}e^{-x^2/4}H_{v}\left(\frac{x}{\sqrt{2}}\right),\quad x
\in \mathbb{R},
\]
in terms of the Hermite polynomial, $H_{\nu}$, of order $\nu$,
whose integral representation is given by
\begin{equation}\label{eq:Hermite}
\mathcal{H}_\nu(z)=\frac{1}{\Gamma(-\nu)}\int_0^\infty
e^{-t^2-2tz}t^{-\nu-1}dt, \quad \text{Re}(\nu)<0.
\end{equation}
 See e.g. Lebedev \cite{lebedev}. The Hermite polynomials satisfy $\mathcal{H}'_\nu(z)=2\nu\mathcal{H}_{\nu-1}(z),
z\in\mathbb{R}$.

\noindent \textbf{Verification of the Conditions in Proposition
\ref{prop:singular-solution}:  }

\noindent (i) Note that
\begin{align}\label{eq:R-F-OU}
R(F(x);b)=\frac{K(x,b)}{\varphi(x)}=\frac{x-b}{\varphi(x)}
\end{align}
and the differentiability of $R(y; b)$ comes from that of
$F(\cdot)$ and $\varphi(\cdot)$.\\ (ii) Since
\begin{equation}
\frac{\partial }{\partial
x}\frac{K(x,b)}{\varphi(x)}=\frac{\varphi(x)-(x-b)\varphi'(x)}{\varphi(x)^2}
\geq 0
\end{equation}
for $x \geq b$, using (\ref{eq:m(x)}), it can be seen that
  $R(y; b)$ is increasing on $y\in
(F(b), \infty)$ to $\infty$.

On the other hand, $(\A-\alpha)(x-b)=p(x)$ for every $x>0$, in
which $p(x)=(1-2\rho x)-\alpha(x-b)$. This linear function $p(x)$
has one positive root at say, $k$.  Then $R(y; \cdot)$ is convex
on $y\in [F(0), F(k))$ and concave on $y\in (F(k), \infty)$.

The facts that $1/\varphi(F^{-1}(y))$ is increasing to $\infty$
and concave on $[F(0),\infty)$ can be similarly shown (by
replacing $K(x,b)$ in (\ref{eq:R-F-OU}) with unity).

\noindent (iii)   We need to verify that (\ref{eq:JS-bstar}) and
(\ref{eq:JS-suf}) hold. In fact, we will see that when $\alpha>0$,
the unique solution of (\ref{eq:JS-bstar}) or the optimal
reflection level satisfies $b^* \in (0,1/(2 \rho))$. This result
is very intuitive, since $x=1/(2 \rho)$ is the mean-reversion
level of $X$.

The functions $\psi(\cdot)$ and $\varphi(\cdot)$ both solve the
differential equation
\begin{equation}\label{eq:diff-eqn}
(1- 2 \rho x) w'(x)+ 2 x fw''(x) -\alpha w(x)=0, \quad
x\in(0,\infty).
\end{equation}
Also, we know from their representation in (\ref{eq:var-phi-psi})
that
\begin{equation}
\psi(x)>0, \,\, \psi'(x)>0; \quad \varphi(x)>0, \,\,
\varphi'(x)<0, \quad x \in (0,\infty).
\end{equation}
Evaluating (\ref{eq:diff-eqn}) at $x=1/(2 \rho)$ for $w=\psi$ and
$w=\varphi$, we obtain
\begin{equation}
2x \psi''\left(\frac{1}{2 \rho}\right)=\alpha \psi\left(\frac{1}{2
\rho}\right) \quad \text{and} \quad 2x \psi''\left(\frac{1}{2
\rho}\right)=\alpha \psi\left(\frac{1}{2 \rho}\right),
\end{equation}
from which it follows that, for $\alpha>0$,
\begin{equation}\label{eq:at-1-2rho}
\frac{\psi''\left(\frac{1}{2 \rho}
\right)}{\varphi''\left(\frac{1}{2 \rho}
\right)}=\frac{\psi\left(\frac{1}{2 \rho}
\right)}{\varphi\left(\frac{1}{2 \rho}  \right)}=F\left(\frac{1}{2
\rho} \right)  \geq F(0)>0,
\end{equation}
where the first inequality follows from the fact that $F$ defined
in (\ref{eq:F}) is increasing. If $\alpha=0$, then
$\psi''(1/(2\rho))=\varphi''(1/(2 \rho))=0$.

First, we will show that (\ref{eq:JS-bstar}) has a unique solution
in $(1, 1/(2 \rho)]$, and show that this solution indeed satisfies
(\ref{eq:JS-suf}). Next, we will show that the same equation does
not have a solution in $(1/(2 \rho),\infty)$. To establish our
first goal, let us collect some information on the behavior of the
functions $\varphi(\cdot)$ and $\psi(\cdot)$ over this interval.
From (\ref{eq:diff-eqn}) with $f=\varphi$, we see that when
\begin{equation}\label{varhidp}
\varphi''(x)>0, \quad x \in \left(0,\frac{1}{2 \rho}\right],
\end{equation}
 since $\varphi'(x)<0$ and $\varphi(x)>0$ over the same
interval. Differentiating (\ref{eq:F}) we obtain
\begin{equation}\label{eq:sec-der}
(2\rho+\alpha) w'(x)-(3-2 \rho  x)w''(x)=  2 x w'''(x),  \quad x
\in (0,\infty),
\end{equation}
for $w=\psi$ or $w=\varphi$. Using (\ref{eq:sec-der}) it can be
seen that
\begin{equation}\label{eq:var-psi-third-sless}
\varphi'''(x)<0, \quad x \in  \left(0,\frac{1}{2 \rho}\right],
\end{equation}
 using the fact that $\varphi'(x)<0$ and
$\varphi''(x)>0$ on the same interval.

After simplifying the expression for $\psi(\cdot)$ in
(\ref{eq:CIR}), we write
\begin{equation}
\psi(x)=H_{-\alpha/\rho}(-\sqrt{\rho
x})-H_{-\alpha/\rho}(\sqrt{\rho x}).
\end{equation}
The second derivative of $\psi(\cdot)$, then can be computed as
\begin{equation}
\psi''(x)=-\frac{\alpha}{4
\sqrt{\rho}}x^{-3/2}\left(H_{-\alpha/\rho-1}(-\sqrt{\rho
x})+H_{-\alpha/\rho-1}(\sqrt{\rho
x})\right)+\frac{\alpha}{4}\left(\frac{\alpha}{\rho}+1\right)x^{-1}
\left(H_{-\alpha/\rho-2}(-\sqrt{\rho
x})+H_{-\alpha/\rho-2}(\sqrt{\rho x})\right),
\end{equation}
from which it follows that
\begin{equation}\label{eq:phi-beh-at-zero}
\lim_{x \rightarrow 0+}\psi''(x)=-\infty.
\end{equation}
With the help of (\ref{eq:sec-der}) with $w=\psi$, observe that
\begin{equation}\label{eq:if-then}
\text{if \, $x \in (0,1/(2 \rho))$, then $\psi'''(x)>0$, if
$\psi''(x)<0$}
\end{equation}
 since $\psi'(x)>0$. It follows from (\ref{eq:at-1-2rho}) and
(\ref{varhidp}) that $\psi''(1/(2 \rho)) > 0$. Now, this fact
together with (\ref{eq:phi-beh-at-zero}) imply that
\begin{equation}\label{eq:psidp-leq}
\psi''(x) \leq 0 \quad \text{for} \,\,\, x\in (0,x_0], \quad
\text{where} \,\, x_0\in \left(0,\frac{1}{2 \rho}\right).
\end{equation}
And it also follows from (\ref{eq:if-then}) that
\begin{equation}\label{eq:psi-dp-pos}
\psi''(x)>0, \quad x \in \left(x_0,\frac{1}{2 \rho}\right).
\end{equation}
At this point, we can state that there exists a solution, $b^* \in
(x_0,1/(2 \rho))$ to (\ref{eq:JS-bstar}) as a result of
(\ref{eq:at-1-2rho}), (\ref{varhidp}), (\ref{eq:psidp-leq}) and
the intermediate value theorem.

Let us prove that at $b^*$, $\psi'''(b^*) \geq 0$. If
$\psi'''(b^*)<0$, then there would exist a point $\tilde{x}_0 \in
(x_0,b^*)$ such that
\begin{equation}\label{eq:psi-t-x-0}
\psi'''(\tilde{x}_0)=0, \quad \psi^{(4)}(\tilde{x}_0)<0,
\end{equation}
in which $\psi^{(4)}$ stands for the fourth derivative of $\psi$.
Differentiating (\ref{eq:sec-der}) we write
\begin{equation}\label{eq:fth-der}
(4 \rho+\alpha)w''(x)-(5-2 \rho x) w'''(x)= 2 x w^{(4)}(x), \quad
x \in (0,\infty).
\end{equation}
Evaluating the left-hand-side of (\ref{eq:fth-der}) (when
$w=\psi$) at $\tilde{x}_0$ we obtain a positive quantity using
(\ref{eq:psi-dp-pos}) and the equality in (\ref{eq:psi-t-x-0}),
whereas the right-hand-side of (\ref{eq:psi-t-x-0}) due to the
inequality in (\ref{eq:psi-t-x-0}), which yields a contradiction.
Since $\psi'''(b^*) \geq 0$, and $\varphi'''(b^*)<0$ by
(\ref{eq:var-psi-third-sless}), $b^*$ satisfies (\ref{eq:JS-suf}).

Let us show that $b^*$ is the only solution to (\ref{eq:JS-bstar})
in $(0,1/(2 \rho))$. Assume there exists another solution to
(\ref{eq:JS-bstar}) in $(0,1/(2 \rho))$, then necessarily there
would be at least one more solution to (\ref{eq:JS-bstar}) in
$(0,1/(2 \rho))$ (This follows from (\ref{eq:at-1-2rho}),
(\ref{varhidp}), (\ref{eq:var-psi-third-sless}), and
(\ref{eq:psidp-leq})). Let us denote the largest three of all of
the solutions by $x_1<x_2<x_3$. It can be easily under this
assumption
\begin{equation}
\psi''(x)>\varphi''(x), \quad x \in (x_1,x_2), \quad \text{and}
\quad \psi'''(x_2)<0.
\end{equation}
But this contradicts the fact we have proved above: since $x_2$ is
a solution to (\ref{eq:JS-bstar}), $\psi'''(x_2) \geq 0$.

It remains to show that (\ref{eq:JS-bstar}) does not have a
solution in $[1/(2 \rho),\infty)$. It is clear from
(\ref{eq:at-1-2rho}) that $x=1/(2 \rho)$ is not a solution of
(\ref{eq:JS-bstar}). Let us assume that for $x>1/(2 \rho)$,
$\psi''(x)=F(0)\varphi''(x)$. Then,
\begin{equation}
-(1-2 \rho x) \psi'(x)+\alpha \psi(x)=-F(0)(1-2 \rho
x)\psi'(x)+F(0) \alpha \varphi(x),
\end{equation}
which implies that
\begin{equation}\label{eq:eq-for-x-geq}
\alpha (F(x)-F(0))=(1-2 \rho x)\left(\frac{\psi'(x)-F(0)
\varphi'(x)}{\varphi(x)}\right).
\end{equation}
Note that the left-hand-side of (\ref{eq:eq-for-x-geq}) is
non-negative because $F$ is increasing. On the other hand the
right-hand-side of (\ref{eq:eq-for-x-geq}) is negative. This
yields a contradiction.

\noindent \textbf{Verification of the Conditions in Proposition
\ref{prop:more-than-local}:  } The only non-trivial condition to
check is whether $v''(x)\le 0$ for $x\in (0, \infty)$. It is clear
from our analysis above that $\psi''(x)-F(0)\varphi''(x)< 0$ on $x
\in (0, b^*)$. The concavity of $v$ follows from
Remark~\ref{rem:verification-remark}-c.

We can determine $\beta^*$ from (\ref{eq:optimal-beta}) and write
down the value function as
\begin{eqnarray*}
v(x)&=& \begin{cases}
                    v_0(x)\triangleq  \beta^* (\psi(x)-F(0)\varphi(x)), &0\leq x \leq b^*, \\
                    v_0(b^*)+x-b^*, &b^*\leq x.
        \end{cases}
\end{eqnarray*}
in which $\psi(x)$ and $\varphi(x)$ are given by
(\ref{eq:SR-fundamental}) with (\ref{eq:Whittaker}). Figure
\ref{fig:SR} illustrates the function $b \rightarrow \beta(b)$,
the value function, $v$, and its derivatives for a special choice
of parameters.
\begin{figure}[h]
\begin{center}
\begin{minipage}{0.45\textwidth}
\centering \includegraphics[scale=0.70]{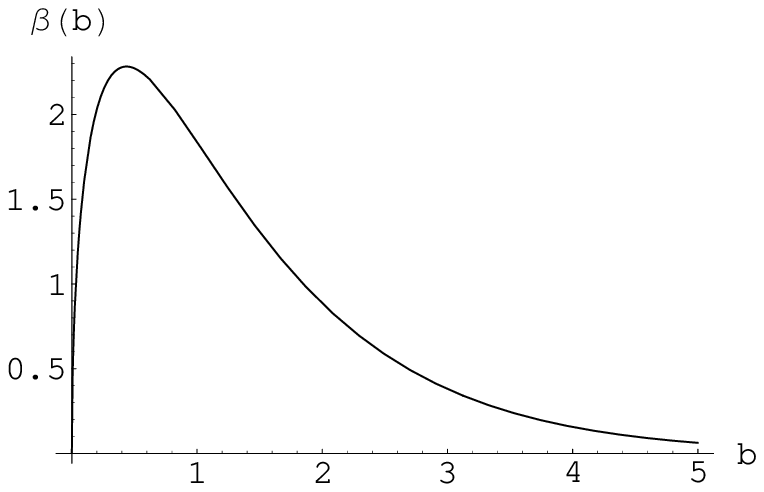} \\ (a)
\end{minipage}
\begin{minipage}{0.45\textwidth}
\centering \includegraphics[scale=0.70]{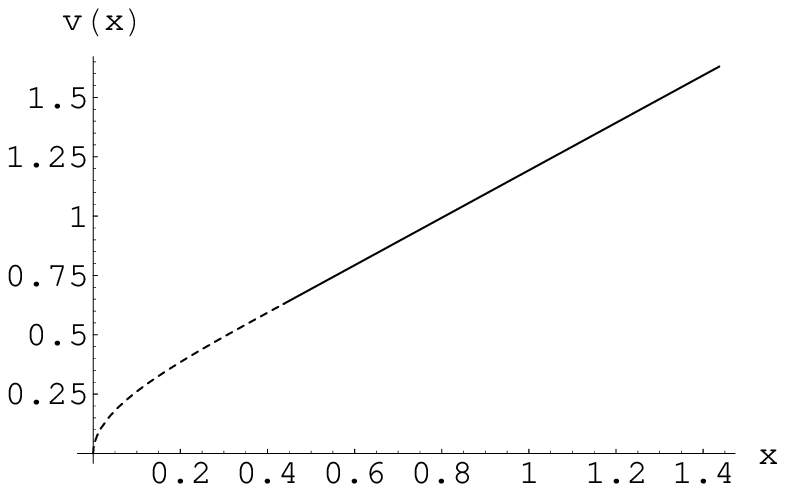} \\ (b)
\end{minipage}
\begin{minipage}{0.45\textwidth}
\centering \includegraphics[scale=0.60]{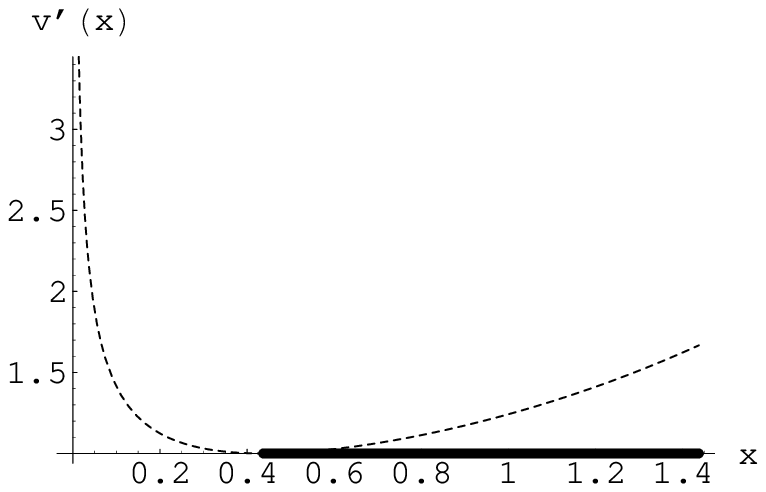} \\ (c)
\end{minipage}
\begin{minipage}{0.45\textwidth}
\centering \includegraphics[scale=0.60]{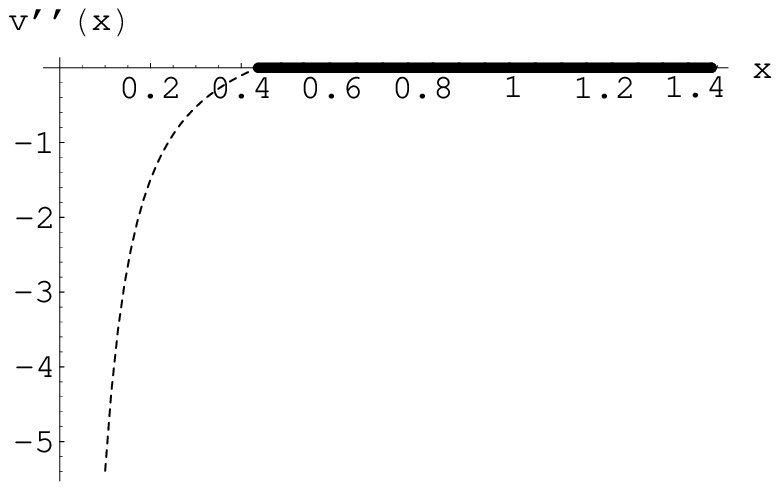} \\ (d)
\end{minipage}
\caption{\small The dividend payout problem with a square root
process with parameters
 $(\rho, \alpha)=(1, 0.2)$:  (a)
 The graph of $b \rightarrow \beta(b)$ (see (\ref{eq:betab})). It attains its maximum at $b^*$.
 (b) The value function $v(x)$. $b^*=0.4370$ and $\beta^*=2.2826$.  (c) The derivative $v'(x)$:
 It is  $v_0'(x)$ on $0\leq x \leq b^*$ and $v'(x)=1$ on $b^*\leq x$.
  The derivative $v_0'(x)$ on $b^*\leq x$ is also shown  to illustrate that $v_0''(b^*)=0$.}
\label{fig:SR}
\end{center}
\end{figure}

\end{example}

\section{Appendix}

\begin{lemma}\label{lem:f-app-smp}
Let us assume that assumption (\ref{eq:f-condition}) holds. Then
for any stopping time $\tau$ of the filtration $(\F_t)_{t \geq 0}$
\begin{equation}
\ME \left[\int_{0}^{\tau} e^{-\alpha s} f(X^0_s)
ds\right]=g(x)-\ME\left[e^{-\alpha \tau} g(X^0_{\tau})\right],
\end{equation}
in which $g$ is defined in (\ref{eq:defn-g}).
\end{lemma}
\begin{proof}
The proof immediately follows from the strong Markov property of
the process $X^0$.
\end{proof}

\begin{lemma}\label{lem:laplace}
For any pair $(l,r) \in (c,d)^2$, let us define
\begin{equation}
v_r(x) \triangleq \ME[e^{-\alpha\tau_r}1_{\{\tau_r<\tau_l\}}],
\quad \text{and} \quad \quad v_l (x)\triangleq
\ME[e^{-\alpha\tau_r}1_{\{\tau_l<\tau_r\}}], \quad x \in [l,r].
\end{equation}
Then
\begin{equation}\label{eq:laplace}
v_r(x)=\frac{\psi(l)\varphi(x)-\psi(x)\varphi(l)}
{\psi(l)\varphi(r)-\psi(r)\varphi(l)}, \quad \text{and} \quad
v_l(x)=\frac{\psi(x)\varphi(r)-\psi(r)\varphi(x)}
{\psi(l)\varphi(r)-\psi(r)\varphi(l)}, \quad x \in [l,r].
\end{equation}
\end{lemma}
\begin{proof}
Both $x \rightarrow v_r(x)$ and $x \rightarrow v_l(x)$, $x \in
(c,d)$ are solutions to $(\mathcal{A}-\alpha)u=0$ with boundary
conditions $v_l(l)=v_r(r)=1$ and $v_l(r)=v_r(l)=0$. Therefore, we
can write them as linear combinations of the homogeneous solutions
of $(\mathcal{A}-\alpha)u=0$, and we get (\ref{eq:laplace}).
\end{proof}

%\noindent\textbf{Acknowledgment} We are grateful to the anonymous
%Associate Editor and the two referees for their comments that
%helped us improve the manuscript.

%\bibliographystyle{plain}
%{\small \bibliography{references}}

\end{document}